\documentclass[leqno,a4paper,oneside]{article}
\usepackage{amsmath,amssymb,amscd,mathrsfs}
\usepackage{graphicx,subfigure,caption,comment,multirow,color}
\usepackage{wrapfig}
\usepackage{lipsum} 
\usepackage{mathptmx,courier,yfonts}
\usepackage[scaled=1.0]{helvet}
\usepackage{longtable}
\begin{document} 

{
\centering
\includegraphics[scale=0.65]{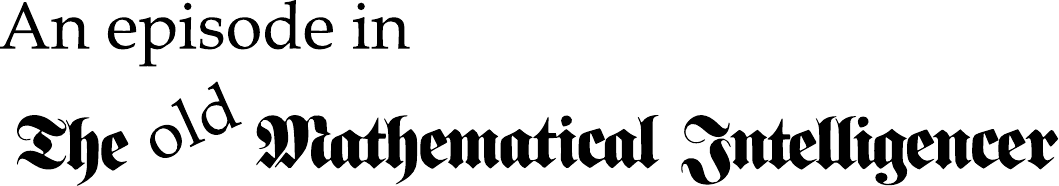}
}

\begin{center}
Christoph Lamm 
\end{center}

\section{Introduction} \label{sec:Introduction}
Recently, in a second-hand bookstore in Germany, I found the bound volume shown on the right.
It contains a complete set of The Mathematical Intelligencer. 
Due to the volume's format and size, this can obviously not mean the currently available mathematical journal of this name!

\begin{wrapfigure}{o}[1.8cm]{0.45\textwidth}
\flushright
\vspace{-0.8cm}
\includegraphics[width=0.4\textwidth]{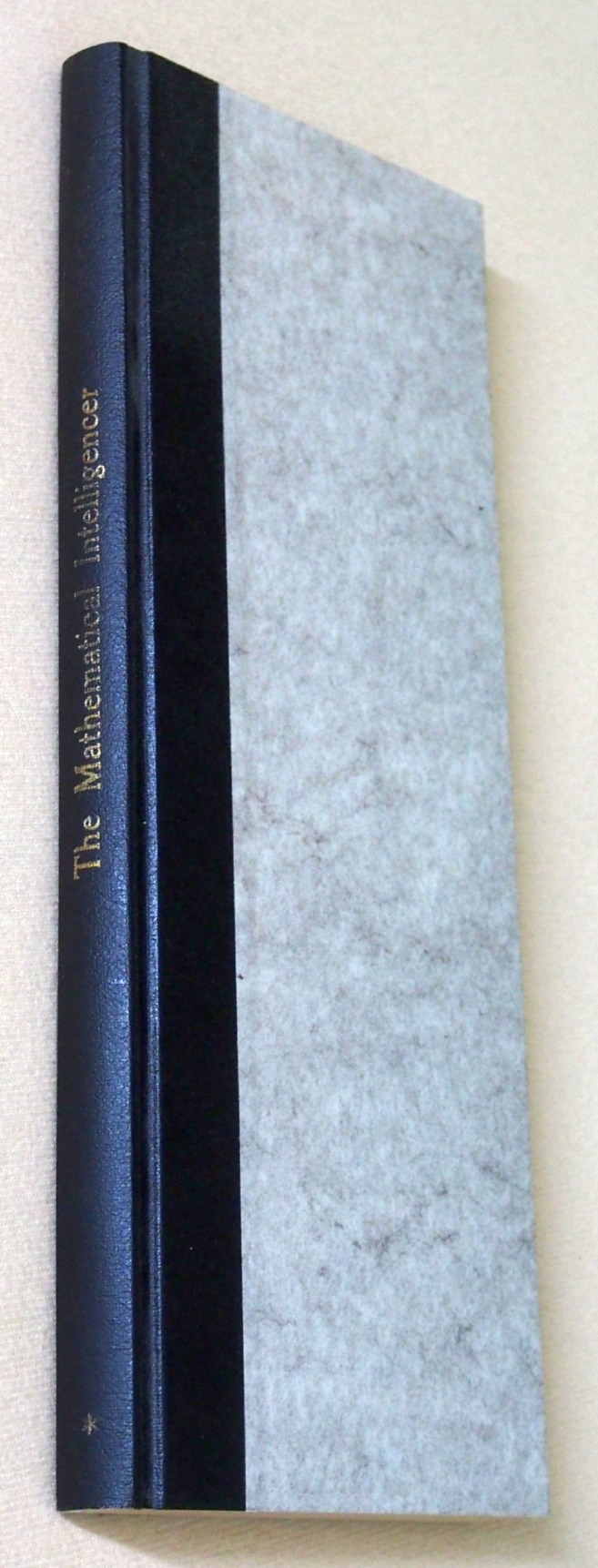}
\label{Bound_volume}
\end{wrapfigure}

Before its Volume 0 appeared in August 1977 (one issue only), followed by the first complete Volume 1 in 1978, 
the Mathematical Intelligencer was a typewritten Springer newsletter, and between 1971 and 1976 in total 12 issues 
appeared (numbered from 0 through 11).
In the issue for the official journal's 40th birthday \cite{Senechal2018}, Marjorie Senechal summarizes:
\begin{quote}
Today’s Mathematical Intelligencer evolved from an occasional typewritten, chatty, foldout newsletter that Springer editors Klaus Peters, Alice Peters, and Walter Kaufmann-Bühler began sending to mathematicians “at home \& abroad” in 1971. As the newsletter grew in size and popularity, Klaus, Alice, and Walter nudged it into a full-fledged magazine, and appointed mathematicians Bruce Chandler and Harold Edwards as co-Editors-in-Chief in 1978. 
\end{quote}

An article in 2008 \cite{Senechal2008} celebrated the 30th birthday and contains more information on the `old' Intelligencer,
for instance on the format and distribution: It was sent as an `accordeon' print to mathematicians world-wide, typically to 7,500.
In particular, Page 8 in \cite{Senechal2008} shows a photo of Alice and Klaus Peters looking at issue 4. We learn from this photo how the 
accordeon print was folded and that one of the contents of issue 4 was a `map of Springer publishing'. This map is the topic of our article.

\clearpage
\newpage
\section{The Realms of Acceptance} \label{sec:Map}
Since the original editorial from issue 4 is now difficult to access, we reproduce it here. 
It provides both background and a taste of the humor in which the map was conceived:

\small
\begin{quote}
WALKS IN MATHEMATICS (Springer Style): The map on pages 3--5 of this Intelligencer was originally drawn by Hilbert in 1921.
Similar to what he did for mathematics in 1900 by his Paris Address, he gave, with this map, an overview of the main endeavors 
of mathematical publishing in the forseeable future. In a famous address to the assembled mathematicians of Göttingen, Berlin,
Heidelberg and Helmstedt (very old university, but rather unknown), F.\,Springer pledged himself and his company to (this one 
of) Hilbert's programs. The first visible consequence was the start of the Grundlehren series in the fall of the same year.

It is perhaps not yet entirely known that the three editors of the M.I. are very much concerned with this map (and, of course,
the appropriate landscape). As keen and professional rangers all three of them know the whole region, but each has his favorite haunts,
where he likes to dwell and to work. People say, in the center of the enchanged island of the Journals, in the middle of a lovely grove of 
evergreen pines there is an old palace; this is the place where Klaus Peters prefers to live and to follow everything that his shining 
eyes can behold from the top of his belfry. He often sees his two colleagues in the parlor of his castle, in the middle of shining armouries
as out of times of old---they have come home, mildly weary from their excursions in the realm. Both have their favorite tracks,
Alice Peters, being a native American, in the Lecture Notes' Republic and, in part, in both the Old Countries, the Kingdom and the Sovereign Duchy.
Walter Kaufmann-Bühler mainly spends his time in the burgeoning Graduate Province and the Twin Counties, but Alice and he keep good 
relations and do not consider excursions into the other's provinces as trespassing, especially since strongholds of either side, 
marked either by A.P. or KB can be found throughout the whole realm. 
A consequence of the whole system is that all three think that all new or important incidents should be (and are) 
discussed among them at a great length and in detail. This, of course, will be slightly more difficult if KB, as he intents to,
will set sail and move his home across the ocean to settle in the New World where the wealth of his two main provinces lies, but regular
visits of (A. and K.) P. to the New World and of KB to the Old one are planned. This is, how the world has been divided and how it is being
cultivated. Sometimes one may be inclined to think it is (or will be) a sheltered park, but there is still a changing number of very
active volcanoes which erupt, laying waste large stretches of land and rendering the soil fertile again.

\flushright{WKB}
\end{quote}
\normalsize

The map is shown on the next page. It is clearly one of the highlights of the old Intelligencer and I studied it closely.
The first part of \cite{Senechal2008} is a conversation between Marjorie Senechal, Alice and Klaus Peters.
(Walter Kaufmann-Bühler had not been present because he had died in 1986, see \cite{MathIntell_1987_4}.)
On the fourth page of this conversation, the map is printed (in a relatively small size) and the figure's caption attributes it 
to Walter Kaufmann-Bühler. However, hidden in the mountaineous parts in the lower middle of the map, I had found the signature 
`Ulrich Felgner'---therefore I suspected that he might have drawn the map. This happened before I looked up other sources, 
as e.g.\,the conversation in \cite{Senechal2008}. Ulrich Felgner confirmed that he drew the map, using the ideas of the three 
Springer editors. In addition to addressing the question of authorship, our article also aims to explain the map.
This will be the subject of the next section.

\enlargethispage{1cm}

\clearpage

\pagenumbering{gobble}

\vspace*{-5cm}\hspace*{-4.7cm}
\includegraphics[angle=270,scale=0.76]{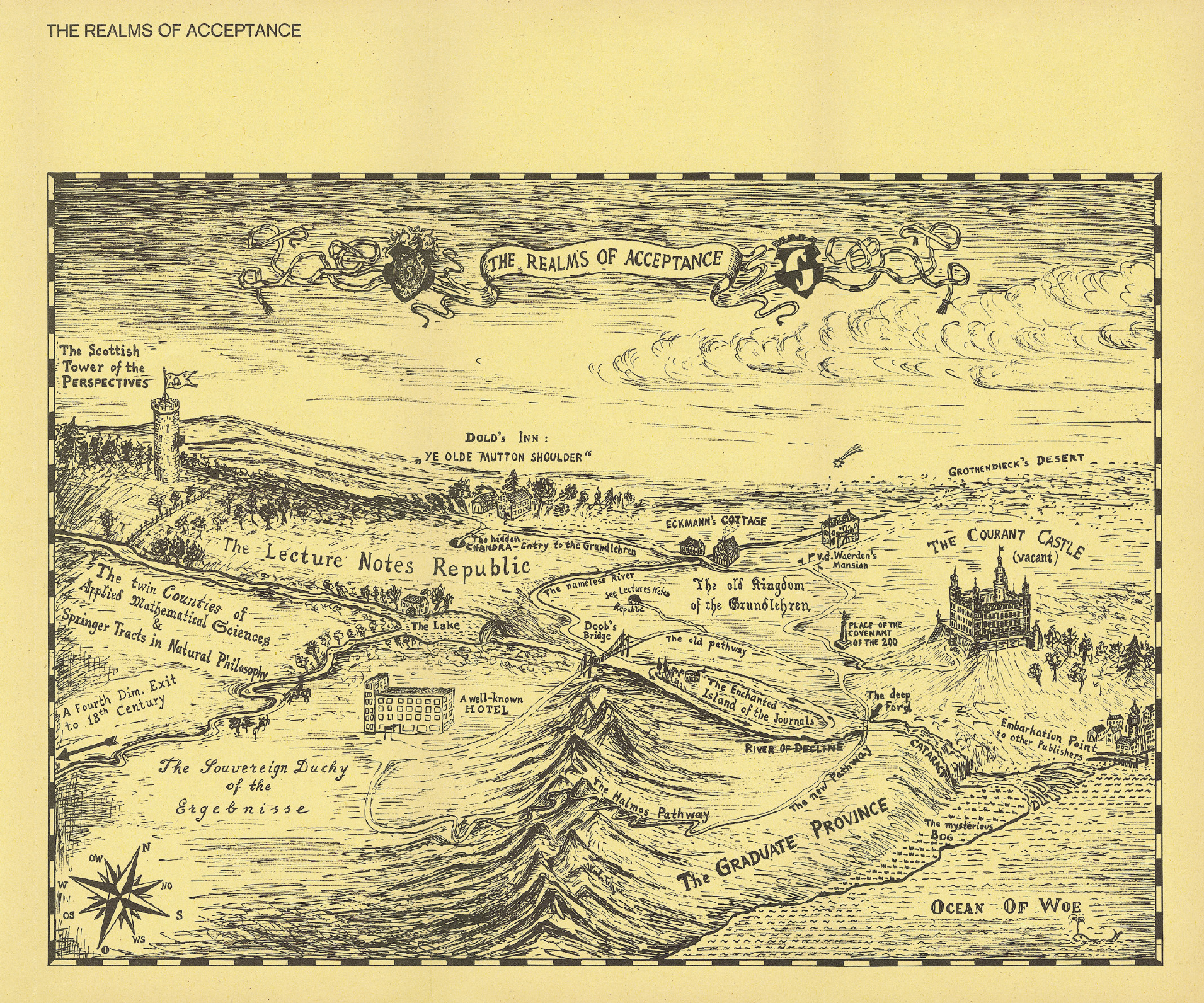}

\section{Annotations} \label{sec:Annotations}
The editorial, reprinted on page \pageref{sec:Map}, already outlines the motivation behind the creation of the map. 
For mathematicians, it is of course evident that David Hilbert was not involved. 
Nearly all of the map's features are now understood and annotated. We begin with the territories.

\enlargethispage{1cm}

\subsubsection*{Territories: Springer's book series and journals}
Annotating the book series and journals is straighforward:
The two oldest book series are the `Grundlehren der mathematischen Wissenschaften' (starting in 1921) and
the `Ergebnisse der Mathematik und ihrer Grenzgebiete' (since 1932), both having several hundreds volumes.
The series `Springer Tracts in Natural Philosophy', more oriented to physics, was started in 1964 and has had 39 volumes.

\begin{wrapfigure}{o}[2cm]{0.45\textwidth}
\flushright
\vspace{-1cm}
\includegraphics[width=0.4\textwidth]{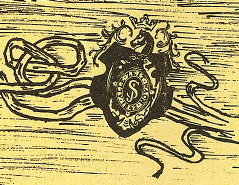}
\label{Logo_left}
\end{wrapfigure}

The first volume of the series `Applied Mathematical Sciences' appeared in 1971, and has now more than 200 volumes.
In the same year, the `Graduate Texts in Mathematics' began, currently with 305 volumes.
The `Perspectives in Mathematical Logic' will have its first book three years in the future (Jon Barwise: Admissible Sets 
and Structures: An Approach to Definability Theory, 1975).
An especially successful series is the `Lecture Notes in Mathematics', beginning in 1964, having now more than 2300 titles.

The mathematical journals are not listed individually in the map. They include for instance the `Mathematische Annalen' 
(since 1920 with Springer). We mention that Klaus Peters initiated the new journal `manuscripta mathematica' in 1969, 
after Karl Stein approached him, see \cite{Peters2009}, p.\,819. 

\pagenumbering{arabic}
\setcounter{page}{4}
				
Less straightforward is the `Fourth Dim.\,Exit to 18th Century' in the left of the map. Readers may submit their interpretations!

\subsubsection*{Castles and cottages}
Among the buildings, `The Courant Castle' is outstanding in its size and its details. It is vacant, though, because Richard Courant 
died in January 1972. The other editors' homes are less glamorous (their castles are of the moderate types `inn', `mansion', `cottage', 
`tower'; and there is also a hotel). We note, that Klaus Peter's home is called an `old palace' in the editorial; it is 
drawn very small in the map, however. 

\begin{itemize}
	\item Richard Courant: We mentioned the vacant castle. Courant remained over 50 years in close contact with Springer Verlag.
	      He was the first editor of the Grund\-lehren series (then in Göttingen). After emigration, he was professor at New York University.
				Its mathematics department is today called the `Courant Institute of Mathematical Sciences'.				
	\item David Hilbert: A well-known hotel---we assume that this is Hilbert's hotel! In a lecture in 1924/25 (`Über das Unendliche'), 
	      Hilbert used a hotel as an illustration for infinite sets. After sketching the finite case, he continues with a hotel
				with infinitely many rooms: `Sobald nun ein neuer Gast hinzukommt, braucht der Wirt nur zu veranlassen, dass jeder der alten Gäste 
				das Zimmer mit der um 1 höheren Nummer bezieht, und es wird für den Neuangekommenen das Zimmer 1 frei.' (\cite{EwaldSieg2013}, p.\,730).
	\item B.L.\,van der Waerden: He was editor of the `Mathematische Annalen' and of the `Grundlehren'. 
	      His influential book `Moderne Algebra' appeared in 1930/31 as `Grundlehren' (Vol.\,33 and 34).
	\item Albrecht Dold: Editor for the `Lecture Notes' (1964--2001), `Ergebnisse' (1968--1982) and `manuscripta mathematica' (1969--1993).
	      He wrote `Grundlehren', Vol. 200 (Lectures on Algebraic Topology). The name `Ye olde mutton shoulder' for Dold’s Inn reflects his 
				reputation as an outstanding cook and generous host \cite{Peters}.
	\item Beno Eckmann: Editor for the `Lecture Notes' and the `Grundlehren'. Professor at the ETH Zürich (1948--1984).
	\item K.\,Chandrasekharan: Professor for number theory at the `Tata Institute of Fundamental Research' (Mumbai) from 1950--1965, and at
	      the ETH Zürich (1965--1988). The map has `The hidden Chandra-Entry to the Grundlehren'. This might refer to his two
				`Grundlehren' volumes 148 and 167, `Introduction to analytic number theory' (1968) and `Arithmetical Functions' (1970).
\end{itemize}

\subsubsection*{The $\Omega$-group}
`The Scottish Tower of the Perspectives' in the upper left side of the map was a mystery to me. However, Ulrich Felgner helped: 
The $\Omega$-group consisted of the logicians Robin Gandy, Hans Hermes, Azriel Lévy, Gert Heinz Müller, 
Gerald Sacks and Dana Scott, who worked on the book series `Perspectives in Mathematical Logic'. 
The name therefore derives from Dana Scott and the name of the book series (and the flag on the tower shows an $\Omega$).

In the map, the lake's mill is intended to allude to Gert Heinz Müller \cite{Felgner}. He organized the `Perspectives' book series and 
also edited the six magnificent volumes of the `$\Omega$-Bibliography of Mathematical Logic', which appeared in 1987 with Springer.

\subsubsection*{Pathways, rivers, lake, bridges, and the Ocean of Woe}
The pathways could be interpreted as the organizational processes in the Springer mathematics department.
The submitted manuscripts might have been assessed in a traditional or more modern way (old and new pathways), or in Halmos style.
When they are rejected, the `River of Decline' takes it down to the `Ocean of Woe'! (Other interpretations are welcome.)

\begin{itemize}
	\item Paul Halmos: Beginning in 1970, he was editor for the `Ergebnisse' and for `Graduate Texts in Mathematics' (\cite{Goetze}, p.\,60).
	      He wrote `How to remember Walter Kaufmann-Bühler' (\cite{MathIntell_1987_4}, p.\,4--10). This article 
				also gives an insight into the collaboration between a publisher’s editor and an academic mathematician.
	\item Lester R.\,Ford Jr.: He studied flow problems and therefore we assume that he was the inspiration for `The deep Ford'.
	      In 1970, Halmos got a Lester Ford Award, an award which honors excellent expository articles in `The American Mathematical Monthly'.
				This award is, however, named after Lester R.\,Ford Sr., and it is today called `Paul R. Halmos -- Lester R. Ford Award'.
	\item Joseph L.\,Doob: Editor for the `Grundlehren'. Open: Did his `Upcrossing inequality' inspire `Doob's bridge'? 
	\item The `Place of the Covenant of the 200' has a triumphal column and celebrates the completion of 200 volumes of the `Grundlehren'
	      (in our context, `covenant' is a humorous term for this achievement).
	\item Alexander Grothendieck: He received the Fields Medal in 1966. Together with Jean Dieudonné, he wrote the `Grundlehren'
	      Vol. 166 (`Éléments de Géométrie Algé\-brique'), appearing in 1971. 
				
				The anecdote explaining `Grothendieck's desert' is the following:
				Desert is a reference to the chaotic state of Grothendieck’s workspace in his apartment. Once, when Klaus Peters and Walter Kaufmann-Bühler
				visited him at home to discuss publishing projects and there was no free table available, Grothendieck simply lifted one side of a large table, 
				letting everything on it slide off onto the floor. The table was now clear, but the floor had become a complete desert \cite{Felgner}.
\end{itemize}

Summarizing, the title `The Realms of Acceptance' covers the manuscript flow in the Springer mathematics departments.
This flow might lead to the acceptance of the manuscript (article, book, lecture note), or to its rejection.
It then could be lost (`The mysterious bog'), or submitted to another publisher (`Embarkation point to other publishers')!

\section{Typewriter typography} \label{sec:typewriter}
The old Mathematical Intelligencer, as well as the Lecture Notes in Mathematics were produced with typewriters.
Everything not available on a typewriter was inserted by hand. Paul Halmos described this `filling in' as follows (\cite{Halmos1985}, p.\,76/77):

\begin{quote}
The meaning of `fill in' is probably no longer known; in those days it was an inevitable part of mathematical typing.
If you were lucky, your typewriter had a plus sign and an equal sign, but Greek letters and inequalities, not to mention
integrals and the other weird symbols that mathematicians are fond of, were unheard of. [\ldots] The result was that
mathematical typing consisted of typing the prose and skipping enough spaces to leave room for letters from other 
alphabetsa and for weird symbols---once that was done you went back over the page and by hand, with a pen, \textsl{filled in}
the subscripts and the rest that you had left out.
\end{quote}

An example from a Lecture Notes volume which appeared in 1971 is shown in the following figure.
Here, the curly P and H, and the direct sum signs are hand-written.

\begin{figure}[hbtp]
\centering
\includegraphics[scale=0.4]{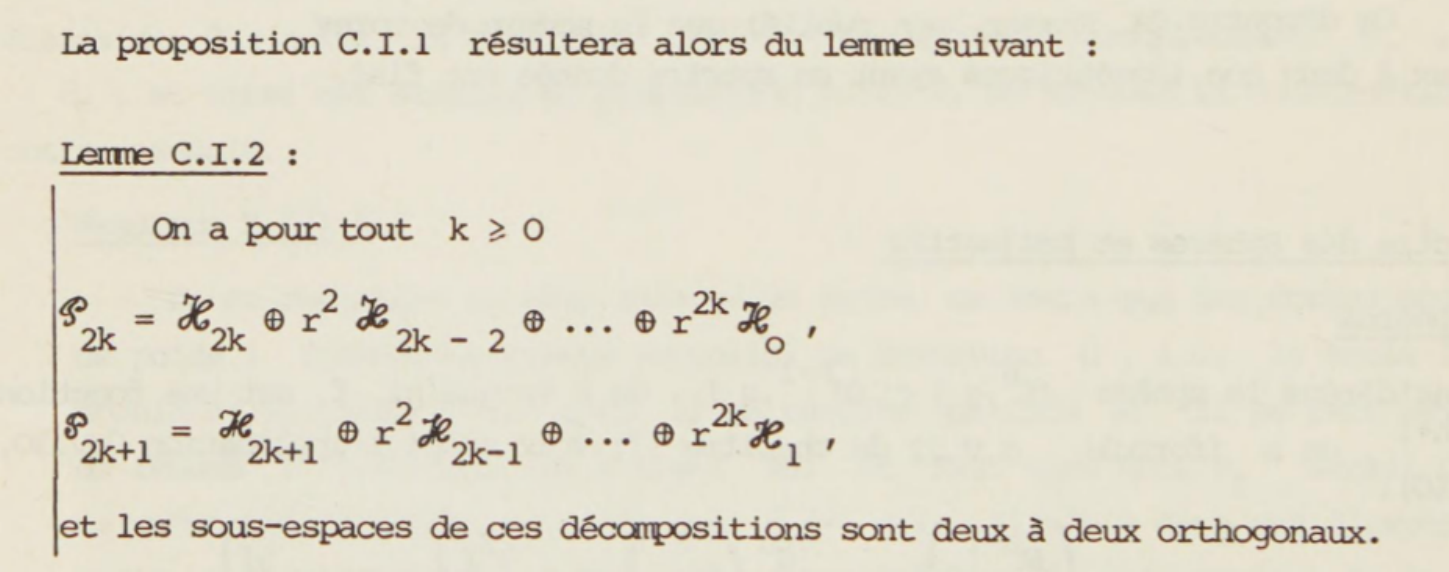}
\label{Berger}
\caption{Part of page 160, in M.\,Berger/P.\,Gauduchon/E.\,Mazet: `Le Spectre d`une Variété Riemannienne', Lecture Notes in Mathematics 194, 1971}
\end{figure}

\newpage
Judging from the year when \cite{Halmos1985} appeared (1985), the use of more flexible typewriters and of course 
of \TeX \,changed that habit already more than 30 years ago.

When I looked up the quotation by Halmos, I remembered that I have his `Finite dimensional vector spaces' in one of 
the first versions from Princeton University Press. Suddenly, I recalled my surprise when I had bought this book, 
also from a second-hand bookshop: The previous owner was K.\,Chandrasekharan!

\begin{figure}[hbtp]
\centering
\includegraphics[scale=0.33]{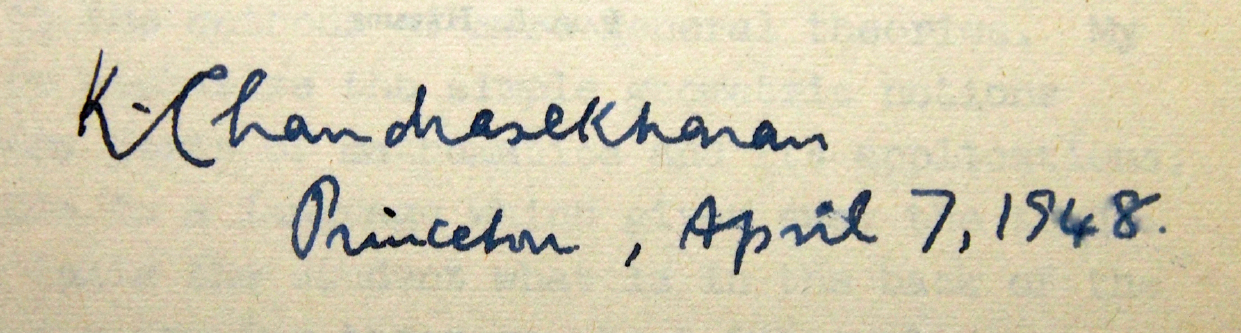}
\label{HalmosChandra}
\caption{Ownership signature by `Chandra' in Halmos' book `Finite dimensional vector spaces'}
\end{figure}

\section*{Acknowledgments and an invitation}
Many thanks to Prof. Ulrich Felgner for sharing memories from 55 years ago and kindly permitting the reprinting of the map, 
and to Alice Peters for the delightful recollection of Albrecht Dold and his culinary talents.

I acknowledge the help of the `Deutsche Museum' \cite{DeutschesMuseum}. They sent me a high quality scan of the map.
The size in the accordeon print is 27x30 cm, whereas the original drawing was done on a paper of approximately
size A3, which is 30x42 cm. 

As already mentioned, for the open issues we ask the readers to contribute their ideas, or also more general reminiscences 
on the old Intelligencer. These will be incorporated in a revised version of this article.

\newpage

\section*{Appendix}
This appendix contains a table of contents of the old Mathematical Intelligencer, to give the reader an impression of all 12 issues.

\small
\begin{longtable}[b]{|p{1.3cm}|p{10cm}|}
\hline
Issue & Content \\ \hline
  0, \hspace{1cm} 1971 & Editorial by K. Peters (pages 1,2); Alfred Brauer on Issai Schur (3,4,5); David Hilbert's voice recording (6); Springer Ad (7,8); Beno Eckmann on Heinz Hopf (11,12,13); review of Constance Reid's `Hilbert' by Carl Ludwig Siegel taken from Zentralblatt (13,14,15); announcement of the next issue (16)  \\ \hline
  1, \hspace{1cm} 1971 or 1972 &  Editorial by Peter Hilton on `Zero growth' for mathematicians (1,2); translation of the George Pólya-Hermann Weyl bet as of 1918 (3,4); chatty remarks on Springer books (5,6); Springer Ad (7,8); 30\% yellow sale (33 books) (11); review of G.T. Kneebone of Hilbert/Paul Bernays `Grundlagen der Mathematik' (Second edition), reprinted from Journal of symbolic logic, 1970 (12,13); K. Peters on the series `Heidelberger Taschenbücher' (14), Springer Ad, list of 24 volumes of Heidelberger Taschenbücher (15); correspondence Hermann Minkowski/Hilbert, Carl Friedrich Gauß \& Hilbert poster (16) \\ \hline
  2, \hspace{1cm} April 1972  & Editorial by W. Kaufmann-Bühler on the history of the `Grundlehren der mathematischen Wissenschaften' (1,2); explanation how Springer's `Mathematical Book Service' works (30\% discount list) and letter by Paul Halmos (3); Tizian: Sacred \& profane love, and remark on an article on Alfred Adler in `The New Yorker' (4,5); Friedrich L. Bauer: `Software engineering', reprinted (6); Springer Ad (7,8); H.\,Götze on Richard Courant who died in 1972 (11,12); IHES, Nicolaas Kuiper new director (13); letter by Egbert Brieskorn on Hilton's editorial, letter by A. Weil (typed without spaces and breaks), letter by E. Wette (15,16) \\ \hline
  3, \hspace{1cm} September 1972 & Editorial by K. Peters (1,2); Zentralblatt's index system explained (3); Oberwolfach: demolition of Lorenzenhof, 5 fotos and a poem by S. MacLane (4,5); Springer Ad (6,7,8); Ms.\,A.\,Dell on Springer's Lecture notes in physics (11); stamp, Goethe (12); mathematical poem by S. Abhyankar, and letter by P. Bernays (13,14); letter by P. Hilton as response to Brieskorn (15,16) \\ \hline
  4, \hspace{1cm} $\sim$Dec. 1972 & Editorial by W. Kaufmann-Bühler on the `Realm of acceptance' map (1,2); the `Realm of acceptance' map, containing the different parts of Springer's maths publications (3,4,5); T.A.A. Broadbent's review of `Analytical inequalities' by D.S. Mitrinovic (6); Springer Ad (7,8); announcement for Noether and Schur poster (11); letters by L.T. Gardner on Titian and by J.D. Zund (12); letter by M. Nagaraj on the Hilton/Brieskorn correspondence (13); story on Venn diagrams (14); nine stamps from Nicaragua (15); quotation/Trotzky library card/Adler quiz (16) \\ \hline
  5, \hspace{1cm} February 1973 & Editorial by W. Kaufmann-Bühler (1,2); Irving Kaplansky: Hilbert's problems (to appear in 1977 as Lecture Notes), excerpt from chapter 7 (3,4,5); review of Saunders Mac Lane's `Categories for the working mathematician' (6); Springer Ad (7,8); `Berufspraxis' of mathematicians, a colloqium series in Bielefeld (11); Alexander Grothendieck's SGA in `Lecture Notes in Mathematics'; letter by Albrecht Dold (12); letter by Ms.\,K.\,Padmavally on the Hilton/Brieskorn correspondence (13,14,15); F. Toth: Figure taken from `Lagerungen in der Ebene auf der Kugel und im Raum'/newspaper article from India/Goethe (16) \\ \hline
  6, \hspace{1cm} April 1973 & Editorial by K. Peters (1,2); Springer's `Mathematical Book Service' (again) explained (3), satirical article on Addison Wesley advertisements in the Notices of the AMS, vol. 20, 1 and 2 (4,5); Gauss's Easter rule (6); Springer Ad (7,8); letter by Armand Borel and A. Weil on the appointment of a professor at the Institute of Advanced Study in Princeton (11); letter by I. Valentine (head of the Springer journal subscription department) (12); letter by B. Crstici (13); chocolate for finding misprints in A. Dold's `Lectures on Algebraic Topology' (14); announcement of 6 new volumes of `Graduate texts in mathematics' (15,16) \\ \hline	
  7, \hspace{1cm} August \hspace{2mm} 1973 & Editorial by A. and K. Peters (1,2) on relevance; the Ergebnisse der Mathematik und ihrer Grenzgebiete report series (3); an illustration of topology as a tree, on the occasion of the book of Boto von Querenburg (autor collective in Bochum) (4); `Book Service' announcements (5); photo of Johann Georg Hagen, taken from Minkowski's letters to Hilbert (6); Springer Ad (7,8); letters by E. Wette and by J.D. Stegeman on the Easter rule (11); letter by J. Staples on the functions of research papers (12,13); letter by H.E. Suchlandt on mathematics as fun (14); letter by S. Swaminathan with a quotation of Alfréd Rényi, letter by P. de Witte on ecological relevance (15); letter by G. Hinrichs with quotations in German (16) \\ \hline
  8, \hspace{1cm} February 1973 & Editorial by W. Kaufmann-Bühler on a future biography of Sofya Kovalevskaya (1,2); review and announcement of the two 2-volume sets on ring theory by C. Faith and by F. Anderson/K. Fuller (3); Springer Verlag at AMS San Francisco Meeting 1974 (4); Springer Ads and announcements (5,6,7,8); description of the series `Applied mathematical sciences' (11,12); letter by Helmut Grunsky, letter by E.L. Poiani, letter by C. Coleman (13); letter by R.A. Rudin, letter by J.S. Golan, a Pierre-Augustin Caron de Beaumarchais quotation used by G. Polya (14); letter by Armand Borel and A. Weil (15); Springer announcement and photo of Felix Klein (16) \\ \hline
  9, \hspace{1cm} August \hspace{2mm} 1973 & Editorial with a longer contribution by P. Hilton, 10 years anniversary of Springer New York (1,2,3); mathematical poem by D.\,Höss (Die Zeit, 26.04.1974), quotation by W.B. Yeats (4); Springer Ad (5,6,7,8); letter by N. Balasubramanian on Venn diagrams (11,12); letter by Heinz Bauer on the Erlanger Programm, letter by F.L. Bauer on the quotation by Beaumarchais and a response by Oskar Perron who also used it in one of his books (13); new Springer series `Undergraduate Texts in Mathematics' (14); Jean Dieudonné on the Erlanger Programm (15,16) \\ \hline
  10, \hspace{1cm} September 1975 & Title page with photo of birthplace of Johannes Müller (1); contents and impressum (2); editorial by K. Peters (3,4); the mathematician and astronom Regiomontanus, Johannes Müller (5); Springer Ad (6); J. Dieudonné on the Weil conjectures (7-21); P. Hilton: book review of Herbert Seifert/William Threlfall: `Lehrbuch der Topologie' (1934) (22,23,24); A. Weil: excerpts in French from his autobiography (25,26,27); Springer Ad (28); article on the US copyright law (29,30) \\ \hline
  11, \hspace{1cm} June 1976 & Title page with picture of Pierre de Fermat (1); contents and impressum (2); editorial (3);  Springer Ad (4); W. Kaufmann-Bühler on Otto Neugebauer: `A History of Ancient Mathematical Astronomy', 3 volumes (5,6); Paulo Ribenboim on the early history of Fermat's last theorem (7-21); corrections to `The Weil conjectures' by Dieudonné (21); Springer Ad (22); W. Kaufmann-Bühler on mathematical typography (23-27); Springer Ad (28); K. Jacobs on James Alexander's horned sphere and Louis Antoine's necklace (29,30)	\\ \hline
\end{longtable}

\normalsize

\newpage

\vspace{1cm}
\noindent
Christoph Lamm \\ \noindent
R\"{u}ckertstr. 3, 65187 Wiesbaden \\ \noindent
Germany \\ \noindent
e-mail: christoph.lamm@web.de

\end{document}